\documentclass[12pt, reqno]{amsart}
\usepackage{pifont}
\usepackage{mathrsfs}
\usepackage{geometry}
\usepackage{titletoc}
\usepackage{stix2}
\usepackage{amscd}
\usepackage{amsmath}
\usepackage{amssymb}
\usepackage{enumitem}
\usepackage{mathtools}
\usepackage[table]{xcolor}
\usepackage[all]{xy}
\usepackage{tikz}
\usepackage{tikz-cd}
\usepackage{indentfirst}
\usepackage{babel}
\usepackage{setspace}
\usepackage[colorlinks,linkcolor=red,anchorcolor=green,citecolor=blue]{hyperref}
\hypersetup{linktocpage = true}
\usepackage{rotating}
\usepackage{ytableau}
\usepackage{longtable}
\usepackage{array}
\usepackage{fancyhdr}
\usepackage{amsthm}

\newcolumntype{M}[1]{>{\centering\arraybackslash}m{#1}}

\newcommand\cA{{\mathcal A}}
\newcommand\cB{{\mathcal B}}

\newcommand\cD{{\mathcal D}}
\newcommand\cE{{\mathcal E}}
\newcommand\cF{{\mathcal F}}

\newcommand\cR{{\mathcal R}}

\newcommand{\id}{\mathrm{Id}}
\newcommand{\rk}{{\rm rk}}

\newcommand\bp{{\bar\partial}}

\theoremstyle{plain}

\newtheorem{thm}{Theorem}[section]
\newtheorem{lemma}[thm]{Lemma}
\newtheorem{prop}[thm]{Proposition}
\newtheorem{cor}[thm]{Corollary}
\newtheorem{defn}[thm]{Definition}

\theoremstyle{definition}
\newtheorem{example}[thm]{Example}
\newtheorem{remark}[thm]{Remark}

\makeatletter
\newcommand{\namedlabel}[2]{\phantomsection\def\@currentlabel{#1}\label{#2}}
\makeatother

\newcommand{\btheorem}{\begin{thm}}
	\newcommand{\etheorem}{\end{thm}}
\newcommand{\bproposition}{\begin{prop}}
	\newcommand{\eproposition}{\end{prop}}
\newcommand{\bdefinition}{\begin{defn}}
	\newcommand{\edefinition}{\end{defn}}
\newcommand{\bcorollary}{\begin{cor}}
	\newcommand{\ecorollary}{\end{cor}}
\newcommand{\bproof}{\begin{proof}}
	\newcommand{\eproof}{\end{proof}}
\newcommand{\bremark}{\begin{remark}}
	\newcommand{\eremark}{\end{remark}}
\newcommand{\bexample}{\begin{example}}
	\newcommand{\eexample}{\end{example}}
\newcommand{\blemma}{\begin{lemma}}
	\newcommand{\elemma}{\end{lemma}}

\newcommand{\p}{\partial}
\renewcommand{\bar}{\overline}

\renewcommand{\phi}{\varphi}

\newcommand{\beq}{\begin{equation}}
\newcommand{\eeq}{\end{equation}}
\newcommand{\be}{\begin{eqnarray*}}
	\newcommand{\ee}{\end{eqnarray*}}
\newcommand{\bd}{\begin{enumerate}}
	\newcommand{\ed}{\end{enumerate}}

\renewcommand{\hat}{\widehat}
\renewcommand{\tilde}{\widetilde}
\newcommand{\qtq}[1]{\quad\mbox{#1}\quad}
\renewcommand{\bp}{\bar{\partial}}
\newcommand{\ts}{\otimes}

\newcommand{\diag}{\operatorname{diag}}
\renewcommand{\>}{\rightarrow}
\newcommand{\C}{{\mathbb C}}
\newcommand{\R}{{\mathbb R}}

\newcommand{\om}{\omega}

\newcommand{\tr}{\mathrm{tr}}

\setlist[itemize]{leftmargin=*}
\setlist[enumerate]{leftmargin=*}

\numberwithin{equation}{section}

\title{Iterative construction of Hermitian-Einstein metrics on stable bundles}
\author{}
\date{}

\author{Jiaxuan Fan}
\author{Zhiyao Xiong}
\author{Xiaokui Yang}
\author{Shing-Tung Yau}

\address{Zhiyao Xiong, Department of Mathematics, Tsinghua University, Beijing, 100084, China}
\email{xiongzy22@mails.tsinghua.edu.cn}
\address{Xiaokui Yang, Department of Mathematics and Yau Mathematical Sciences Center, Tsinghua University, Beijing, 100084, China}
\email{xkyang@mail.tsinghua.edu.cn}
\address{Jiaxuan Fan, Qiuzhen College, Tsinghua University, Beijing, 100084, China}
\email{fanjiaxu21@mails.tsinghua.edu.cn}
\address{Shing-Tung Yau, Yau Mathematical Sciences Center and Qiuzhen College, Tsinghua University, Beijing, 100084, China}
\email{styau@mail.tsinghua.edu.cn}

\begin{document}
	
	\begin{abstract}
		Let $E$ be a stable  holomorphic vector bundle over a compact K\"ahler (or  Gauduchon) manifold $(M,\omega_g)$. We show that
		for any real number $\mu>0$ and any initial Hermitian metric $h_0$ on $E$, there exists a unique iteration sequence $\{h_m\}$ satisfying
		$$
		\Lambda_{\omega_g}\left(\sqrt{-1}R^{h_{m+1}}\right)
		=(\lambda_E-\mu)h_{m+1}+\mu h_m,
		$$
		 and  $\{h_m\}$ converges smoothly to a Hermitian--Einstein metric $h_\infty$ on $E$ satisfying
		$$
		\Lambda_{\omega_g}\left(\sqrt{-1}R^{h_{\infty}}\right)
		=\lambda_Eh_\infty,
		$$
		where $\lambda_E\in \R$ is the stability constant. A key feature of this proof is that it is independent of Donaldson's variational 
		framework and applies to non-K\"ahler manifolds.
	\end{abstract}
	\maketitle
	
	{\small{    \begin{spacing}{1.1} \tableofcontents %
			\dottedcontents{section}[1.8cm]{}{3em}{5pt} %
\end{spacing} }} 
\vspace*{-1\baselineskip} 
	
	\section{Introduction}
	This paper constitutes a continuation of the program initiated in \cite{WYY26+, FWYY26+, WYY26b+,XYY26+, XYY26b+}.  We refer the reader to these works for a comprehensive account of the historical development of this subject over the past six decades.\\

	Our main result establishes the assertion announced in \cite[Remark~1.8, Item~$4$]{WYY26+} and provides a constructive proof of the Donaldson-Uhlenbeck-Yau theorem (\cite{Don85,UY86,Don87}, see also \cite{NS65} and \cite{LY86}):
	
	\btheorem\label{main theorem}
	Let $E$ be a stable holomorphic vector  bundle over a compact K\"ahler (or Gauduchon) manifold $(M,\omega_g)$.  Then, for any real number $\mu>0$ and any initial Hermitian metric $h_0$ on $E$,  there exists a unique iteration sequence $\{h_m\}$ satisfying
	\beq 
	\Lambda_{\omega_g}\left(\sqrt{-1}R^{h_{m+1}}\right)
	=(\lambda_E-\mu)h_{m+1}+\mu h_m,
	\eeq 
	 and  $\{h_m\}$ converges smoothly to a Hermitian-Einstein metric $h_\infty$ on $E$ satisfying
	\beq 
	\Lambda_{\omega_g}\left(\sqrt{-1}R^{h_{\infty}}\right)
	=\lambda_Eh_\infty,
	\eeq 
	where $\lambda_E\in \R$ is the stability constant
	\begin{equation}
	\lambda_E=\frac{2\pi n\int_M c^{\mathrm{BC}}_1(E)\wedge\omega_g^{n-1}}{\rk(E)\int_M\om_g^n}.
	\label{lambda definition}
	\end{equation}
	\etheorem

	\bremark The sign of $\lambda_E$
	  leads to three natural variants of the iteration.
	\bd
	\item If $\lambda_E>0$, one can choose $\mu=\lambda_E$, reducing the iteration to 
	\beq \Lambda_{\omega_g}\left(\sqrt{-1}R^{h_{m+1}}\right)
	=\lambda_E  h_m.\eeq 
	\item If  $\lambda_E=0$, then the iteration equation becomes
	\beq \Lambda_{\omega_g}\left(\sqrt{-1}R^{h_{m+1}}\right)
	=-\mu h_{m+1}+\mu h_m,\eeq 
	and $\{h_m\}$ converges smoothly to a Hermitian-Einstein metric  $h_\infty$ on $E$ satisfying 
	\beq 
	\Lambda_{\omega_g}\left(\sqrt{-1}R^{h_{\infty}}\right)
	=0.
	\eeq 
	\item If $\lambda_E<0$, then $\lambda_{E^*}>0$ for the dual bundle $E^*$.  One can iterate on $E^*$:
	\beq \Lambda_{\omega_g}\left(\sqrt{-1}R^{h^*_{m+1}}\right)
	=\lambda_{E^*}  h^*_m.\eeq 
	The resulting limit $h^*_\infty$ is a smooth Hermitian-Einstein metric on $E^*$,
	and its dual induces a smooth Hermitian-Einstein metric on $E$.
	\ed
\eremark	
\noindent In particular,	one gets the following special case: 
	\bcorollary\label{special}
	Let $E$ be a stable holomorphic vector bundle over a compact K\"ahler manifold $(M,\omega_g)$ with positive stability constant $\lambda_E$.  Then for any initial Hermitian metric $h_0$ on $E$,   there exists a unique iteration sequence $\{h_m\}$ satisfying
	\beq 
	\Lambda_{\omega_g}\left(\sqrt{-1}R^{h_{m+1}}\right)
	=\lambda_E h_m,
	\eeq 
	and  $\{h_m\}$ converges smoothly to a Hermitian-Einstein metric $h_\infty$ on $E$ satisfying
	$$
	\Lambda_{\omega_g}\left(\sqrt{-1}R^{h_{\infty}}\right)
	=\lambda_Eh_\infty.
	$$
	\ecorollary
	
	\bremark The main results hold for stable Higgs bundles over any compact Gauduchon manifold $(M, \omega_g)$, as established in \cite{FWYY26+, WYY26b+, XYY26b+}.
	\eremark
	
		\bremark For iterative constructions of K\"ahler-Einstein metrics on Fano manifolds, see 
		\cite{Rub07, Rub08, Kel09, DR19}; for the constant scalar curvature case, we refer 
		to \cite{Zha25} and the references therein.
	\eremark
	
	\bremark  During the preparation of this manuscript, we learned that Corollary \ref{special} 
	had been established in \cite{CSZ26+}.
	\eremark
	
\noindent\textbf{Acknowledgements}. The third-named author would like to thank  Bing-Long Chen, Jixiang Fu,  Kefeng Liu and Song Sun  for inspiring  discussions.

\vskip  2\baselineskip
\section{A priori estimates}\label{background section}
Let $E$ be a holomorphic vector bundle of rank $r$ over a compact Gauduchon manifold $(M,\om_g)$, i.e. $\p\bp\omega_g^{n-1}=0$.
For a Hermitian metric $h$ on $E$, we write
\begin{equation}
S^h:=\Lambda_{\omega_g}\left(\sqrt{-1}R^h\right)
\in\Gamma(M,E^*\ts\bar E^*),
\qquad
K^h:=S^h\cdot h^{-1}\in\Gamma(M,E^*\ts E).
\label{S and K}
\end{equation}
Suppose that $\{h_0,h_1,\cdots\}$ is a sequence of Hermitian metrics on $E$ satisfying
\beq 
\Lambda_{\omega_g}\left(\sqrt{-1}R^{h_{m+1}}\right)
=(\lambda_E-\mu)h_{m+1}+\mu h_m. \label{iteration tensor equation}
\eeq 
If we write
\begin{equation}
H_m=h_m\cdot h_0^{-1},
\qquad
Q_m=H_{m-1}\cdot H_m^{-1}=h_{m-1}\cdot h_m^{-1},
\label{H Q definition}
\end{equation}
for $m\geq 1$, then the equation \eqref{iteration tensor equation} is equivalent to
\begin{equation}
K^{h_m}=(\lambda_E-\mu)\id_E+\mu Q_m.
\label{iteration endomorphism equation}
\end{equation}
Clearly,  the  endomorphism $Q_m\in \Gamma(M,E^*\ts E)$ is positive and self-adjoint with respect to both $h_{m-1}$ and $h_m$.

\blemma\label{adjacent comparison lemma}
There exists a constant $C=C(h_0,h_1,\lambda_E ,\mu)$ such that
\begin{equation}
C^{-1}h_m\leq h_{m-1}\leq Ch_m
\label{adjacent comparison estimate}
\end{equation}
for all $m\geq 1$. 
In particular,
\begin{equation}
C^{-1}\id_E\leq Q_m\leq C\id_E,
\label{Qbounds}
\end{equation} with respect to $h_m$.
Moreover, one has
\begin{equation}
\|K^{h_m}\|_{C^0(M,\om_g,h_m)}\leq C,
\label{K bounds}
\end{equation}
and 
\beq \lim_m \|Q_{m+1}-Q_m\|_{L^2(M,\om_g,h_m)}=0. \label{dm tends zero}\eeq 
\elemma

\begin{proof}  It is clear that there is a constant $C_1=C_1(h_0,h_1)$ such that 
	\eqref{adjacent comparison estimate} holds for $m=1$ with $C=C_1$. We show by induction that  	\eqref{adjacent comparison estimate} holds for all $m\geq 1$  with this constant $C$.  Suppose that \eqref{adjacent comparison estimate} holds for $m$ with constant $C$.   By this induction hypothesis, one has
	\beq 
	S^{h_m}-(\lambda_E-\mu) h_m=\mu h_{m-1}
	\leq C\mu h_m
	=S^{Ch_{m+1}}- (\lambda_E-\mu) \left(Ch_{m+1}\right).
	\eeq 
	By \cite[Theorem~2.3]{WYY26b+}, one yields
	\beq h_m\leq C h_{m+1}.\eeq   Similarly,
	\beq 
	S^{h_{m+1}}-(\lambda_E-\mu)  h_{m+1}=\mu h_m
	\leq C\mu h_{m-1}
	=S^{Ch_m}-(\lambda_E-\mu)  (Ch_m),
	\eeq 
	and so
	$h_{m+1}\leq C h_m$.  Thus the estimate \eqref{adjacent comparison estimate} holds	for  $m+1$. The equation
	\eqref{iteration endomorphism equation} gives
	\beq 
	\|K^{h_m}\|_{C^0(M,\om_g,h_m)}
	\leq \sqrt r\left(|\lambda_E-\mu|+\mu C\right).
	\eeq 
	
	Moreover, if we use the Hermitian-Yang-Mills functional
	\beq  \cE(h)=\int_M\tr_E\left((K^h-\lambda_E\id_E)^2\right)\omega_g^n,\eeq 
	and write $\cE_m:=\cE(h_m)$, then it is known that 
	\beq \cE_m-\cE_{m+1}\geq \|K^{h_m}-K^{h_{m+1}}\|_{L^2(M,\om_g,h_m)}^2 \label{energy dissipation estimate}\eeq 
	for all $m\geq 0$. Indeed, since both $K^{h_m}$ and $K^{h_{m+1}}$ are self-adjoint with respect to $h_m$, 	a direct calculation shows
	\beq
	\cE_m-\cE_{m+1}
	=\|K^{h_m}-K^{h_{m+1}}\|_{L^2(M,\om_g,h_m)}^2 +2\left(K^{h_m}-K^{h_{m+1}},
	K^{h_{m+1}}-\lambda_E\id_E\right)_{g,h_m}.
	\label{energy polarization}
	\eeq 
	By \cite[Proposition~3.6]{WYY26+}, one has 
	\beq 
	K^{h_{m+1}}-K^{h_m}
	=\sqrt{-1}\Lambda_{\omega_g}\bar\partial
	\left(\partial^{h_m}(Q_{m+1}^{-1})\cdot Q_{m+1}\right)
	=-\sqrt{-1}\Lambda_{\omega_g}\bar\partial
	\left(Q_{m+1}^{-1}\cdot\partial^{h_m}Q_{m+1}\right).
	\eeq  Since equation \eqref{iteration endomorphism equation} has the form
	\beq 
	K^{h_{m+1}}-\lambda_E\id_E=\mu(Q_{m+1}-\id_E),
	\eeq
	we obtain
	\begin{eqnarray}
	&&\left(K^{h_m}-K^{h_{m+1}},
	K^{h_{m+1}}-\lambda_E\id_E\right)_{g,h_m}
	\nonumber\\
	&=&\mu\int_M
	\left\langle
	\sqrt{-1}\Lambda_{\omega_g}\bar\partial
	\left(Q_{m+1}^{-1}\cdot\partial^{h_m}Q_{m+1}\right),
	Q_{m+1}-\id_E
	\right\rangle_{h_m}\omega_g^n
	\nonumber\\
	&=&\mu\int_M
	\left\langle
	Q_{m+1}^{-1}\cdot\partial^{h_m}Q_{m+1},
	\partial^{h_m}Q_{m+1}
	\right\rangle_{g,h_m}\omega_g^n\nonumber\\
	&&+\mu n\sqrt{-1}\int_M
	\tr_E\left(Q_{m+1}^{-1}\cdot\partial^{h_m}Q_{m+1}\cdot(Q_{m+1}-\id_E)\right)
	\wedge\bar\partial\omega_g^{n-1}.
	\label{positive cross term with torsion}
	\end{eqnarray}
	For the torsion term, we compute
	\begin{eqnarray}
	&&\tr_E\left(Q_{m+1}^{-1}\cdot\partial^{h_m}Q_{m+1}\cdot(Q_{m+1}-\id_E)\right)
	\nonumber\\
	&=&
	\tr_E\left(Q_{m+1}^{-1}\cdot\partial^{h_m}Q_{m+1}\cdot Q_{m+1}\right)
	-\tr_E\left(Q_{m+1}^{-1}\cdot\partial^{h_m}Q_{m+1}\right)
	\nonumber\\
	&=&
	\tr_E\left(\partial^{h_m}Q_{m+1}\right)
	-\partial\log\det Q_{m+1}
	\nonumber\\
	&=&
	\partial\left(\tr_E Q_{m+1}-\log\det Q_{m+1}\right).
	\label{trace torsion identity}
	\end{eqnarray}
	Here we use $\tr_E(\partial^{h_m}Q_{m+1})=\partial\tr_E Q_{m+1}$ and  Jacobi's formula
	\[
	\partial\log\det Q_{m+1}
	=\tr_E\left(Q_{m+1}^{-1}\cdot\partial^{h_m}Q_{m+1}\right).
	\]
	Thus Stokes' theorem and the Gauduchon condition
	$\partial\bar\partial\omega_g^{n-1}=0$  yield
	\[
	\int_M
	\partial\left(\tr_E Q_{m+1}-\log\det Q_{m+1}\right)
	\wedge\bar\partial\omega_g^{n-1}=0.
	\]
	Combining this with \eqref{positive cross term with torsion}, we obtain
	\beq
	\left(K^{h_m}-K^{h_{m+1}},
	K^{h_{m+1}}-\lambda_E\id_E\right)_{g,h_m}
	=\mu\int_M
	\left|Q_{m+1}^{-1/2}\cdot\partial^{h_m}Q_{m+1}
	\right|_{g,h_m}^2\omega_g^n\geq0.
	\label{positive cross term}
	\eeq
	This proves \eqref{energy dissipation estimate}. If the  equality in \eqref{energy dissipation estimate} holds, then $\partial^{h_m}Q_{m+1}=0$ and so  $K^{h_m}=K^{h_{m+1}}$.  By \eqref{iteration endomorphism equation}, one has
	\beq 
	K^{h_{m+1}}-K^{h_m}=\mu(Q_{m+1}-Q_m).
	\eeq
	In particular, for all $m\geq 1$,
	\beq \|Q_{m+1}-Q_m\|^2_{L^2(M,\om_g,h_m)}=\mu^{-2} \|K^{h_{m+1}}-K^{h_m}\|^2_{L^2(M,\om_g,h_m)}\leq \mu^{-2}\left( \cE_m-\cE_{m+1}\right).\eeq 
	Since the non-negative sequence $\{\cE_m\}$ is non-increasing,  one obtains \eqref{dm tends zero} and this completes the proof.
\end{proof}

\noindent 
For $m\geq 0$, we set 
\begin{equation}
F_m:=\int_M\log\det H_m\,\omega_g^n,
\label{Fm definition}
\end{equation}
and 
\begin{equation}
c_m:=\exp\left(-\frac{F_m}{r\int_M\om_g^n}\right),
\qquad
\hat h_m:=c_mh_m,
\qquad
\hat H_m:=c_mH_m.
\label{normalization definition}
\end{equation}
It is clear that 
\begin{equation}
\int_M\log\det\hat H_m\,\omega_g^n=0.
\label{normalized determinant mean}
\end{equation}
\blemma\label{Fm monotonicity lemma}
\bd
\item The sequence $\{F_m\}$ is non-decreasing, and $F_m\geq0$.
\item There exists a positive constant $C=C(h_0,h_1,\lambda_E,\mu)$ such that for $m\geq 1$
\beq
C^{-1}c_{m-1}\leq c_m\leq c_{m-1}.
\eeq
\ed
\elemma

\begin{proof}
	Since $Q_m=H_{m-1}H_m^{-1}$,
	\begin{equation}
	F_m-F_{m-1}=-\int_M\log\det Q_m\,\omega_g^n.
	\label{F difference}
	\end{equation}
	Taking the trace of \eqref{iteration endomorphism equation} and integrating, we obtain
	\begin{equation}
	\int_M\tr_EQ_m\,\omega_g^n=r\int_M\om_g^n.
	\label{Q trace average}
	\end{equation}
	Since $\log$ is concave, Jensen's inequality together with \eqref{Q trace average}
	gives
	\begin{eqnarray}
	\int_M\log\det Q_m\,\omega_g^n
	&\leq& r\int_M\log\frac{\tr_EQ_m}{r}\,\omega_g^n\nonumber\\
	&\leq&\left( r\int_M\om_g^n\right)\cdot 
	\log\left(
	\frac{1}{r\int_M\om_g^n}\int_M\tr_EQ_m\,\omega_g^n
	\right)=0.
	\end{eqnarray}
	Hence $F_m\geq F_{m-1}$, and $F_0=0$. This implies (1). On the other hand, 
	\beq \frac{c_m}{c_{m-1}}=\exp\left(-\frac{F_m-F_{m-1}}{r\int_M\om_g^n}\right)=\exp\left(\frac{\int_M\log\det Q_m\,\omega_g^n}{r\int_M\om_g^n}\right).\eeq 
	By \eqref{Qbounds}, we obtain $(2)$. 
\end{proof}

\blemma\label{matrix remainder lemma}
For positive-definite $r\times r$ Hermitian matrices $A,B$  and $0<\sigma\leq1$, define
\begin{align}
\cD_\sigma(A,B)
&:=\sigma\tr\left(AB^{\sigma-1}-B^\sigma\right)
-\tr\left(A^\sigma-B^\sigma\right),
\label{D sigma definition}\\
\cR_\sigma(A,B)
&:=\tr\left(A^{1+\sigma}B^{-1}+B^\sigma-A^\sigma-AB^{\sigma-1}\right).
\label{R sigma definition}
\end{align}
Then 
\begin{equation}
0\leq\cD_\sigma(A,B)
\leq(1-\sigma)\cR_\sigma(A,B).
\label{matrix remainder inequality}
\end{equation}
\elemma

\begin{proof}
	Choose unitary diagonalizations
	\[
	A=U\diag(a_1,\cdots,a_r)U^*,
	\qquad
	B=V\diag(b_1,\cdots,b_r)V^*,
	\]
	and put $W=U^*V=(w_{ik})$.  Direct calculation gives
	\begin{align}
	\cD_\sigma(A,B)
	={}&\sum_{i,k}|w_{ik}|^2b_k^\sigma
	\left[
	\sigma\left(\frac{a_i}{b_k}-1\right)
	-\left(\left(\frac{a_i}{b_k}\right)^\sigma-1\right)
	\right],
	\label{D double sum}\\
	\cR_\sigma(A,B)
	={}&\sum_{i,k}|w_{ik}|^2b_k^\sigma
	\left(\frac{a_i}{b_k}-1\right)
	\left(\left(\frac{a_i}{b_k}\right)^\sigma-1\right).
	\label{R double sum}
	\end{align}
	It is therefore enough to prove, for every $x>0$,
	\begin{equation}
	0\leq\sigma(x-1)-(x^\sigma-1)
	\leq(1-\sigma)(x-1)(x^\sigma-1).
	\label{scalar remainder inequality}
	\end{equation}
	The first inequality follows from the concavity of $x^\sigma$.  For the second one, one has
	\beq 
	(1-\sigma)(x-1)(x^\sigma-1)
	-\left[\sigma(x-1)-(x^\sigma-1)\right]
	=x^\sigma\left((1-\sigma)x+\sigma-x^{1-\sigma}\right)\geq0,
	\eeq
	because $x^{1-\sigma}\leq(1-\sigma)x+\sigma$. 
\end{proof}

\noindent
We use the normalization:
\begin{equation}
\Lambda_m:=\sup_M\lambda_{\max}(H_m),
\qquad
\tilde H_m:=\Lambda_m^{-1}H_m.
\label{Lambda tilde definition}
\end{equation}

\btheorem\label{record index thm}
Let $m_i\to\infty$ satisfy
\begin{equation}
\Lambda_{m_i}=\max_{0\leq j\leq m_i}\Lambda_j.
\label{record index condition}
\end{equation}
Then for every $0<\sigma\leq1$,
\begin{equation}
\limsup_{i\to\infty}
\int_M\tr_E\left(\left(Q_{m_i}-\id_E\right)\cdot\tilde H_{m_i}^\sigma\right)\omega_g^n
\leq0.
\label{record index conclusion}
\end{equation}
Moreover,
\begin{equation}
\limsup_{i\to\infty}
\left[
\int_M\tr_E\left(K^{h_{m_i}}\cdot\tilde H_{m_i}^\sigma\right)\omega_g^n
-\lambda_E\int_M\tr_E\left(\tilde H_{m_i}^\sigma\right)\omega_g^n
\right]\leq0.
\label{record curvature estimate}
\end{equation}
\etheorem

\begin{proof}
	Fix $\sigma\in(0,1]$.   For $m\geq1$, we set
	\beq
	T_m:=\int_M\tr_E\left((Q_m-\id_E)\cdot\tilde H_m^\sigma\right)\omega_g^n.
	\eeq
	Then \eqref{record index conclusion} is equivalent to
	\beq
	\limsup_{i\to\infty}T_{m_i}\leq0.
	\label{T record index conclusion}
	\eeq
	Define
	\begin{align}
	\cA_m&:=\int_M\tr_E\left(\left(Q_m-\id_E\right)\cdot H_m^\sigma\right)\omega_g^n,
	\qquad m\geq1,\label{A m definition}\\
	\cB_m&:=\int_M\tr_E\left(H_m^\sigma\right)\omega_g^n,
	\qquad m\geq0,\label{B m definition}\\
	\varepsilon_m&:=\int_M\tr_E\left(\left(Q_{m+1}-Q_m\right)\cdot H_m^\sigma\right)\omega_g^n,
	\qquad m\geq1.\label{epsilon m definition}
	\end{align}
	For $m\geq1$,  we set
	\begin{equation}
	\cD_m:=\sigma\cA_m-(\cB_{m-1}-\cB_m).
	\label{D m definition}
	\end{equation}
	Since both $H_m$ and $H_{m-1}$ are self-adjoint with respect to $h_0$, we  apply Lemma~\ref{matrix remainder lemma}  with $A=H_{m-1}$ and $B=H_m$. This yields
	\begin{equation}
	0\leq\cD_m\leq(1-\sigma)\cR_m,\qquad m\geq1,
	\label{D R m inequality}
	\end{equation}
	where
	\[
	\cR_m:=\int_M\tr_E\left(
	H_{m-1}^{1+\sigma}\cdot H_m^{-1}+H_m^\sigma
	-H_{m-1}^\sigma-H_{m-1}\cdot H_m^{\sigma-1}
	\right)\omega_g^n\geq0.
	\]
	For $m\geq2$, a direct calculation gives
	\begin{equation}
	\cR_m=\cA_{m-1}-\cA_m+\varepsilon_{m-1}.
	\label{R m telescoping identity}
	\end{equation}
	Hence, for $m\geq2$,
	\begin{equation}
	\cA_{m-1}\geq\cA_m-|\varepsilon_{m-1}|.
	\label{backward propagation inequality}
	\end{equation}
	We claim that there exists a constant
	$C_1=C_1(\sigma,\omega_g,h_0,h_1,\lambda_E,\mu)>0$ such that, whenever
	$2\leq p\leq q$ and
	$\Lambda_q=\max_{0\leq j\leq q}\Lambda_j$,
	\begin{equation}
	\sum_{m=p}^q\cA_m
	\leq C_1\Lambda_q^\sigma
	\left(1+\sum_{j=p-1}^{q-1}
	\|Q_{j+1}-Q_j\|_{L^2(M,\om_g,h_j)}\right).
	\label{A sum upper bound}
	\end{equation}
	Indeed,  for all $j\leq q$, $H_j\leq\Lambda_q\id_E$, and
	\beq 
	0\leq H_j^\sigma\leq\Lambda_q^\sigma\cdot\id_E
	\eeq 
	with respect to $h_j$.
	Together with the estimate \eqref{Qbounds} in Lemma~\ref{adjacent comparison lemma}, we deduce that there exists a constant $C_2=C_2(\omega_g,h_0,h_1,\lambda_E,\mu)>0$ such that for all $j\leq q$,
	\begin{equation}
	|\cA_j|+\cB_j\leq C_2\Lambda_q^\sigma,
	\qtq{and}
	|\varepsilon_j|
	\leq C_2\|Q_{j+1}-Q_j\|_{L^2(M,\om_g,h_j)}\Lambda_q^\sigma.
	\label{A epsilon estimate}
	\end{equation}
	By \eqref{D m definition}, \eqref{D R m inequality} and \eqref{R m telescoping identity}, we have
	\begin{eqnarray}
	\sigma\sum_{m=p}^q\cA_m
	&=&\cB_{p-1}-\cB_q+\sum_{m=p}^q\cD_m\nonumber\\
	&\leq&\cB_{p-1}
	+(1-\sigma)\left(\cA_{p-1}-\cA_q+
	\sum_{m=p}^q\varepsilon_{m-1}\right).
	\label{summed remainder inequality}
	\end{eqnarray}
	Substituting \eqref{A epsilon estimate} into \eqref{summed remainder inequality} and dividing both sides by $\sigma$, this yields a constant $C_1=C_1(\sigma,\omega_g,h_0,h_1,\lambda_E,\mu)>0$ such that \eqref{A sum upper bound} holds.\\

	Suppose for contradiction that \eqref{T record index conclusion} fails.
	Since $T_m=\Lambda_m^{-\sigma}\cA_m$, the failure of
	\eqref{T record index conclusion} implies that, after possibly passing to a
	subsequence, there exists $\delta>0$ such that
	\begin{equation}
	\cA_{m_i}\geq\delta\Lambda_{m_i}^\sigma.
	\label{record contradiction assumption}
	\end{equation}
	Define \beq C_3:=C_1\cdot\left(1+\frac{\delta}{2C_2}\right).\eeq  Choose $\rho>0$
	sufficiently small so that
	\begin{equation}
	\frac{\delta s_\rho}{2}>C_3\qtq{where}
	s_\rho:=\left\lfloor\frac{\delta}{2C_2\rho}\right\rfloor.
	\label{s rho definition}
	\end{equation}
	By \eqref{dm tends zero} in
	Lemma~\ref{adjacent comparison lemma}, we can choose $N_\rho\geq1$ such
	that for all $m\geq N_\rho$,
	\beq
	\|Q_{m+1}-Q_m\|_{L^2(M,\om_g,h_m)}\leq\rho.\label{delta Q rho}
	\eeq
Fix a sufficiently large $i$ so that
	$m_i-s_\rho\geq N_\rho$.  Choose an integer $m$ with $$m_i-s_\rho+1\leq m\leq m_i.$$  
	For every $j\in\{m,\cdots,m_i-1\}$, since
	$\Lambda_{m_i}=\max_{0\leq \ell\leq m_i}\Lambda_\ell$, 
	\eqref{A epsilon estimate} and \eqref{delta Q rho} give
	\beq
	|\varepsilon_j|
	\leq C_2\|Q_{j+1}-Q_j\|_{L^2(M,\om_g,h_j)}\Lambda_{m_i}^\sigma
	\leq C_2\rho\Lambda_{m_i}^\sigma.
	\eeq
	Combining this with \eqref{backward propagation inequality}, we obtain
	\beq
	\cA_m
	\geq \cA_{m_i}-\sum_{j=m}^{m_i-1}|\varepsilon_j|
	\geq \delta\Lambda_{m_i}^\sigma
	-s_\rho C_2\rho\,\Lambda_{m_i}^\sigma
	\geq\frac{\delta}{2}\Lambda_{m_i}^\sigma
	\eeq
	where the last inequality follows from
	$s_\rho C_2\rho\leq\delta/2$.
	Since $m$ is arbitrary in this range,
	\begin{equation}
	\sum_{m=m_i-s_\rho+1}^{m_i}\cA_m
	\geq\frac{\delta s_\rho}{2}\Lambda_{m_i}^\sigma.
	\label{A sum lower bound}
	\end{equation}
	On the other hand, applying \eqref{A sum upper bound} with
	$p=m_i-s_\rho+1$ and $q=m_i$ gives
	\begin{equation}
	\sum_{m=m_i-s_\rho+1}^{m_i}\cA_m
	\leq C_1\Lambda_{m_i}^\sigma(1+s_\rho\rho)
	\leq C_3\Lambda_{m_i}^\sigma,
	\label{A sum final upper}
	\end{equation}
	where in the second inequality we use the definition of $C_3$. 
	Comparing \eqref{A sum lower bound} and
	\eqref{A sum final upper} gives \beq \delta s_\rho/2\leq C_3\eeq  and this contradicts 
	the choice of $\rho$.  Therefore
	\eqref{T record index conclusion} holds, which implies \eqref{record index conclusion}.\\
	
	Finally, \eqref{iteration endomorphism equation} gives
	\beq 
	K^{h_{m_i}}-\lambda_E\id_E=\mu(Q_{m_i}-\id_E),
\eeq 
and this implies
	\be 
	&&\int_M\tr_E\left(K^{h_{m_i}}\cdot\tilde H_{m_i}^\sigma\right)\omega_g^n
	-\lambda_E\int_M\tr_E\left(\tilde H_{m_i}^\sigma\right)\omega_g^n\nonumber\\
	&=&\int_M\tr_E\left((K^{h_{m_i}}-\lambda_E\id_E)\cdot
	\tilde H_{m_i}^\sigma\right)\omega_g^n
	=\mu\int_M\tr_E\left((Q_{m_i}-\id_E)\cdot
	\tilde H_{m_i}^\sigma\right)\omega_g^n
	=\mu T_{m_i}.
	\ee
	Thus \eqref{record curvature estimate} follows from \eqref{T record index conclusion}.
\end{proof}

\vskip  1\baselineskip

\section{The proof of Theorem \ref{main theorem}}

This section is devoted to the proof of Theorem \ref{main theorem}. We work on a compact Gauduchon manifold  $(M,\omega_g)$ and begin by establishing a few technical preliminaries.
We recall the notations required for the iteration:
\begin{equation}\label{iteration eq recall}
\Lambda_{\omega_g}\!\left(\sqrt{-1}\,R^{h_{m+1}}\right)
= (\lambda_E-\mu)\,h_{m+1}+\mu\,h_m.
\end{equation}
Set
\begin{gather}
H_m := h_m\cdot h_0^{-1},
\qquad
Q_m := H_{m-1}\cdot H_m^{-1},
\label{C0 H Q recall}
\\[6pt]
c_m :=
\exp\!\left(
-\frac{1}{\rk(E)\int_M\omega_g^{\,n}}
\int_M\log\det H_m\;\omega_g^{\,n}
\right),\!
\qquad
\hat{h}_m := c_m h_m,
\qquad
\hat{H}_m := c_m H_m,
\label{C0 normalization recall}
\\[6pt]
\Lambda_m := \sup_M\lambda_{\max}(H_m),
\qquad
\tilde{H}_m := \Lambda_m^{-1}H_m,
\qquad
\hat{\Lambda}_m := \sup_M\lambda_{\max}(\hat{H}_m)
= c_m\Lambda_m.
\label{C0 Lambda recall}
\end{gather}
\blemma\label{determinant oscillation lemma}
There exists a constant
$C=C(\om_g,h_0,h_1,\lambda_E,\mu)$ such that
\begin{equation}
\operatorname{osc}_M\log\det H_m\leq C,
\label{determinant oscillation}
\end{equation}
and
\begin{equation}
e^{-C}\leq\det\hat H_m\leq e^{C}
\label{normalized determinant pointwise}
\end{equation}
for all $m\geq0$.  
\elemma

\begin{proof}
	A straightforward computation  shows
	\beq
	\Delta_{\C}\log\det H_m
	=\tr_EK^{h_0}-\tr_EK^{h_m}.
	\eeq
	By Lemma~\ref{adjacent comparison lemma},  there is a
	constant $C_1=C_1(\omega_g, h_0,h_1,\lambda_E,\mu)$ such that
	\beq
	\|K^{h_m}\|_{C^0(M,\om_g,h_m)}\leq C_1
	\eeq
	for every $m$.  Since $K^{h_m}$ is self-adjoint with respect to $h_m$,
	\beq
	\left|\tr_EK^{h_m}\right|
	\leq\sqrt r\cdot|K^{h_m}|_{h_m}
	\leq C_1\sqrt r.
	\eeq
	Therefore, there exists a constant $C_2=C_2(\om_g,h_0,h_1,\lambda_E,\mu)$ such that 
	\beq
	\left|\Delta_{\C}\log\det H_m\right|\leq \left|\tr_EK^{h_0}\right|+\left|\tr_E K^{h_m}\right|\leq C_2.
	\eeq 
	A standard Green function estimate for the complex Laplacian $\Delta_\C$ on the compact Gauduchon manifold $(M,\omega_g)$ yields 
	\beq
	\operatorname{osc}_M \log\det H_m \leq C_3
	\eeq 
	for some constant $C_3=C_3(\om_g,h_0,h_1,\lambda_E,\mu)$. This establishes the estimate in \eqref{determinant oscillation}.  Since $\hat H_m = c_m H_m$,  one obtains
	\beq 
	\operatorname{osc}_M \log\det \hat H_m=
	\operatorname{osc}_M \log\det H_m \leq C_3.
	\eeq 
	It follows from  \eqref{normalized determinant mean} that 
	\beq
	-C_3\leq \log\det \hat H_m\leq C_3,
	\eeq
	which proves \eqref{normalized determinant pointwise}.
\end{proof}

\blemma\label{lem priori estimates for weak compactness}
For all
$m\geq0$ and $\sigma\in(0,1]$, one yields
\beq
\left| \p^{h_0}\tilde H_m^\sigma \cdot \tilde H_m^{-\sigma/2} \right|_{g,h_0}^2
\leq  \frac{1}{\sigma}\Delta_\C\left(\tr_E \tilde H_m^\sigma\right)
+\tr_E\left(\left(K^{h_m}-K^{h_0}\right)\cdot \tilde H_m^\sigma\right). \label{estimateA}
\eeq
Moreover, there exists a constant $C=C(\omega_g, h_0,h_1,\lambda_E,\mu)$ such that 
\beq
\left|\p^{h_0}\tilde H_m^\sigma\cdot\tilde H_m^{-\sigma/2}\right|_{g,h_0}^2
\leq\frac{1}{\sigma}\Delta_\C\left(\tr_E\tilde H_m^\sigma\right)
+C\tr_E\tilde H_m^\sigma,
\label{tilde H inequality}
\eeq

\beq
\left\|\tilde H_m^\sigma\right\|_{W^{1,2}(M,\omega_g,h_0)}\leq C, 
\label{tilde H W2p estimate}
\eeq
and
\beq
\left\|\tilde H_m\right\|_{L^2(M,\omega_g,h_0)}\geq C^{-1}.
\label{non-trivial condition}
\eeq
\elemma

\begin{proof}
	The estimate in \eqref{estimateA} follows from \cite[Lemma~4.2]{XYY26b+}.
	By Lemma~\ref{adjacent comparison lemma}, there exists a constant $C_1=C_1(\omega_g, h_0,h_1,\lambda_E ,\mu)$ such that for all $m\geq 0$
	\beq
	\|K^{h_m}\|_{C^0(M,\om_g,h_m)}\leq C_1.\label{curvature input}
	\eeq
	As in the proof of \cite[Propositions~4.3~and~4.4]{XYY26b+}, the remaining estimates now follow from \eqref{estimateA} and \eqref{curvature input}.
\end{proof}
\noindent Let $W^{1,2}(M,E^*\ts E)$ be the space of $W^{1,2}$-sections of $E^*\ts E$ with respect to $\omega_g$ and $h_0$.
 Recall that (\cite{UY86}) an element $\pi\in W^{1,2}(M,E^*\ts E)$ is called a weakly holomorphic projection if the following identities
\begin{equation}
\pi^*=\pi=\pi^2,
\qquad
(\id_E-\pi)\circ\bar\partial\pi=0
\label{weak projection identities}
\end{equation}
hold almost everywhere on $M$,  where the adjoint is taken with respect to $h_0$.

\blemma\label{power projection extraction}
Let $m_i\to\infty$ satisfy the condition
\eqref{record index condition} and assume that
$\hat\Lambda_{m_i}\to\infty$.  Set
\begin{equation}
\tilde H_i:=\tilde H_{m_i}=\Lambda_{m_i}^{-1}H_{m_i}.
\label{tilde H i proposition definition}
\end{equation}
There exist a subsequence of $\{\tilde{H}_i\}$, which we still denote by 
$\{\tilde{H}_i\}$, and a sequence $\{\sigma_k\} \subset (0,1]$ decreasing to $0$ such that:
\bd
\item $\tilde H_i\> \tilde H_\infty$ in the weak $W^{1,2}(M,E^*\ts E)$ sense for some nonzero $\tilde H_\infty$ in $W^{1,2}(M,E^*\ts E)$.
Moreover, for each fixed $\sigma_j$, $\tilde H_i^{\sigma_j}\> \tilde H^{\sigma_j}_\infty$ in the weak $W^{1,2}(M,E^*\ts E)$ sense.
\item $\tilde H_\infty^{\sigma_j}\> \tilde H$ in the weak $W^{1,2}(M,E^*\ts E)$ sense.
\item $\pi:=\id_E-\tilde H$ is a weakly holomorphic projection of $E$.
\item Let $\cF\subset E$ be the coherent subsheaf associated with $\pi$. Then $0<\rk(\cF)<\rk(E)$.
\ed
\elemma

\bproof  As in the proof of \cite[Theorem~5.2]{XYY26b+}, one can prove (1)--(3) using Lemma~\ref{lem priori estimates for weak compactness}. By Lemma~\ref{Fm monotonicity lemma}, the sequence $\{F_m\}_{m\geq0}$ is non-decreasing and $F_m\geq0$.  Hence
\beq
0<c_m=\exp\left(-\frac{F_m}{r\int_M\om_g^n}\right)\leq1,
\qquad
\hat\Lambda_m=c_m\Lambda_m\leq\Lambda_m.
\eeq
The assumption $\hat\Lambda_{m_i}\>\infty$ implies
$\Lambda_{m_i}\>\infty$.  The normalization gives 
\begin{equation}
\tilde H_i=\hat\Lambda_{m_i}^{-1}\hat H_{m_i}.
\label{tilde normalized identity}
\end{equation}
Thus, by \eqref{normalized determinant pointwise},
\begin{equation}
0<\det\tilde H_i
=\hat\Lambda_{m_i}^{-r}\det\hat H_{m_i}
\leq e^{C}\hat\Lambda_{m_i}^{-r}\>0
\label{tilde determinant vanishing}
\end{equation}
uniformly on $M$.   Since $\tilde H_i\>\tilde H_\infty$ almost everywhere
along a further subsequence, it follows that
$\det\tilde H_\infty=0$ almost everywhere.  Applying
\cite[Lemma~5.1]{XYY26b+} to $\pi$, we obtain a coherent subsheaf
$\cF\subset E$ and an analytic subset $\Sigma\subset M$.  On
$M\setminus\Sigma$, the projections $\pi$ and
$\tilde H=\id_E-\pi$ have constant ranks.  Then $\tilde H_\infty\neq0$ implies that the projection $\tilde H$ has positive rank on $M\setminus\Sigma$, while $\det\tilde H_\infty=0$ ensures that this rank is strictly smaller than $r$.  This completes the proof.
\eproof

\bproposition\label{iteration proposition}
Let $E$ be a stable holomorphic vector bundle over a compact Gauduchon manifold $(M,\om_g)$.
Fix $\mu>0$ and an initial metric $h_0$ on $E$.
Then there exists a unique sequence of Hermitian metrics $\{h_m\}_{m\geq 0}$ on $E$ satisfying
\begin{equation}
S^{h_{m+1}}=\left(\lambda_E-\mu\right) h_{m+1}+\mu h_m.
\end{equation}
With the notations introduced above, there exists a constant $\tilde C_0=\tilde C_0(\om_g,h_0,\lambda_E,\mu)$ such that for all $m\geq 0$, 
\beq
\tilde C_0^{-1}h_0\leq \hat h_m\leq \tilde C_0h_0.\label{normalized C0 estimate}
\eeq 
Moreover, for every integer $k\geq 1$ there exists a constant $\tilde C_k=\tilde C_k(k,\om_g,h_0,\lambda_E,\mu)$ such that for all $m\geq 0$,
\begin{equation}
\|\hat h_m\|_{C^k(M,\om_g,h_0)}\leq \tilde C_k.
\label{normalized Ck estimate}
\end{equation}

\eproposition

\begin{proof}
	Since $E$ is stable, 
	\beq \lambda_E^-=\lambda_E>\lambda_E-\mu.\eeq  
	By \cite[Theorem~1.1]{WYY26b+}, for any $m\geq 0$, there exists a unique  solution $h_{m+1}$ to the equation \eqref{iteration tensor equation}.  Suppose that 
	\beq  \limsup_m \hat\Lambda_m=\infty. \eeq   Hence, there exists a sequence $\{m_i\}\>\infty$ such that
	\begin{equation}
	\hat \Lambda_{m_i}=\max_{0\leq j\leq m_i}\hat \Lambda_j, \qtq{and} \lim_i  	\hat \Lambda_{m_i}=\infty.
	\label{record choice C0}
	\end{equation}
	Since $c_m$ is positive and  non-increasing (Lemma \ref{Fm monotonicity lemma}) and $\hat\Lambda_{m}=c_{m}\Lambda_{m}$, one has 
	\begin{equation}
\Lambda_{m_i}=\max_{0\leq j\leq m_i} \Lambda_j.
	\label{hat Lambda record blowup}
	\end{equation}
	Set
	\begin{equation}
	\tilde H_i:=\tilde H_{m_i}=\Lambda_{m_i}^{-1}H_{m_i}.
	\label{tilde H i}
	\end{equation}   Hence, there exist a sequence
	 $\{\sigma_k\}\subset (0,1]$ decreasing to $0$,  and 
	$\tilde H_\infty,\tilde H\in W^{1,2}(M,E^*\ts E)$  satisfying the properties in Lemma \ref{power projection extraction}.  In particular, $\pi:=\id_E-\tilde H$ is a weakly holomorphic projection.  This projection
	determines a coherent subsheaf $\cF\subset E$ satisfying
	$0<\rk(\cF)<r$.	Since $0<\tilde H_i\leq\id_E$, one has
	$\tilde H_i^{-\sigma_k}\geq\id_E$, and 
	\beq 
	\left|\partial^{h_0}\tilde H_i^{\sigma_k}\right|_{g,h_0}^2
	\leq
	\left|\partial^{h_0}\tilde H_i^{\sigma_k}\cdot
	\tilde H_i^{-\sigma_k/2}\right|_{g,h_0}^2 .
	\eeq
	Then by the estimate in \eqref{estimateA}, we deduce that
	\beq 
	\int_M\left|\partial^{h_0}\tilde H_i^{\sigma_k}\right|_{g,h_0}^2
	\omega_g^n
	+\int_M\tr_E\left(K^{h_0}\cdot\tilde H_i^{\sigma_k}\right)\omega_g^n
	\leq
	\int_M\tr_E\left(K^{h_{m_i}}\cdot\tilde H_i^{\sigma_k}\right)\omega_g^n.
	\label{integrated degree power inequality}
	\eeq
	 The weak
	$W^{1,2}$ convergence of
	$\tilde H_i^{\sigma_k}\rightarrow \tilde H_\infty^{\sigma_k}$ implies
	\begin{equation}
	\int_M\left|\partial^{h_0}\tilde H_\infty^{\sigma_k}\right|_{g,h_0}^2
	\omega_g^n
	\leq
	\limsup_{i\to\infty}
	\int_M\left|\partial^{h_0}\tilde H_i^{\sigma_k}\right|_{g,h_0}^2
	\omega_g^n,
	\label{Dirichlet lower semicontinuity}
	\end{equation}
 and
	\begin{equation}
	\lim_{i\to\infty}\int_M\tr_E\left(K^{h_0}\cdot\tilde H_i^{\sigma_k}\right)\omega_g^n
	=
	\int_M\tr_E\left(K^{h_0}\cdot \tilde H_\infty^{\sigma_k}\right)\omega_g^n.
	\label{K h0 zeroth order convergence}
	\end{equation}
	Moreover, by Theorem~\ref{record index thm},
	\begin{equation}
	\limsup_{i\to\infty}
	\int_M\tr_E\left(K^{h_{m_i}}\cdot\tilde H_i^{\sigma_k}\right)\omega_g^n
	\leq
	\lambda_E\int_M\tr_E\left(\tilde H_\infty^{\sigma_k}\right)\omega_g^n.
	\label{record curvature fixed exponent}
	\end{equation}
	Combining \eqref{integrated degree power inequality},
	\eqref{Dirichlet lower semicontinuity},
	\eqref{K h0 zeroth order convergence}, and
	\eqref{record curvature fixed exponent}, we obtain
	\beq
	\int_M\left|\partial^{h_0}\tilde H_\infty^{\sigma_k}\right|_{g,h_0}^2
	\omega_g^n
	+\int_M\tr_E\left(K^{h_0}\cdot \tilde H_\infty^{\sigma_k}\right)\omega_g^n
	\leq
	\lambda_E\int_M\tr_E\left(\tilde H_\infty^{\sigma_k}\right)\omega_g^n.
	\label{fixed exponent degree inequality}
	\eeq
 Since $\tilde H_\infty^{\sigma_k}\to \tilde H$ weakly in $W^{1,2}(M,E^*\ts E)$, we obtain
	\begin{eqnarray}
	\int_M\left|\partial^{h_0}\tilde H\right|_{g,h_0}^2\omega_g^n
	&\leq& \limsup_{k\to\infty}\int_M\left|\partial^{h_0}\tilde H_\infty^{\sigma_k}\right|_{g,h_0}^2\omega_g^n\nonumber\\
	&\leq& \limsup_{k\to\infty}\left( \lambda_E\int_M\tr_E\left(\tilde H_\infty^{\sigma_k}\right)\omega_g^n - \int_M\tr_E\left(K^{h_0}\cdot \tilde H_\infty^{\sigma_k}\right)\omega_g^n\right)\nonumber\\
	&=&  \lambda_E\int_M \rk\left(\tilde H\right)\om_g^n -\int_M \tr_E\left(K^{h_0}\cdot \tilde H\right)\omega_g^n.
	\label{projection derivative degree estimate}
	\end{eqnarray}
	By  similar arguments as in the proof of
	\cite[Theorem~1.1]{XYY26b+}, we have
	\[
	\deg_{\omega_g}(\cF)
	=\int_M\tr_E(K^{h_0}\cdot\pi)\,\omega_g^n
	-\int_M|\partial^{h_0}\pi|_{g,h_0}^2\,\omega_g^n.
	\]
	Combining this identity with \eqref{projection derivative degree estimate}, one has
	\begin{eqnarray}
	\deg_{\omega_g}(\cF)
	&\geq&
	\int_M\tr_E(K^{h_0}\cdot\pi)\,\omega_g^n
	+\int_M \tr_E\left(K^{h_0}\cdot \tilde H\right)\omega_g^n-\lambda_E\int_M \rk\left(\tilde H\right)\om_g^n
	\nonumber\\
	&=&\deg_{\omega_g}(E)
	-\lambda_E\cdot \left(r-\rk(\cF)\right)\cdot \int_M\omega_g^n\nonumber\\
	&=&\lambda_E\cdot \rk(\cF)\cdot \int_M\omega_g^n.
	\label{stability contradiction inequality}
	\end{eqnarray}
This contradicts the stability of $E$, since $0 < \operatorname{rk}(\mathcal{F}) < r$.  Therefore,
\beq  \limsup_m \hat\Lambda_m\leq C_1, \eeq  
where	 $C_1=C_1(\om_g,h_0,\lambda_E,\mu)$.  Moreover, by
\eqref{normalized determinant pointwise}, there is a constant
$C_2=C_2(\om_g,h_0,\lambda_E,\mu)$ such that
\beq
e^{-C_2}\leq \det\hat H_m\leq e^{C_2}\label{det hat H bound}
\eeq
for all $m\geq 0$. 
	Hence, there exists a constant $\tilde C_0=\tilde C_0(\om_g,h_0,\lambda_E,\mu)$ such that the estimate in \eqref{normalized C0 estimate} holds for all $m\geq 0$.\\
	
	For higher-order estimates, by Lemma~\ref{Fm monotonicity lemma}, for all $m\geq 0$,
	\beq
	C_3^{-1}\leq\frac{c_{m+1}}{c_{m}}\leq1,\label{cm bounds}
	\eeq
	for a constant $C_3=C_3(\omega_g, h_0,\lambda_E,\mu)$.  Multiplying 
	\eqref{iteration tensor equation} by $c_{m+1}$  yields
	\begin{equation}
	S^{\hat h_{m+1}}=(\lambda_E-\mu)\hat h_{m+1}+\mu\frac{c_{m+1}}{c_{m}}\hat h_{m}.
	\label{normalized iteration equation}
	\end{equation}
	The estimates \eqref{normalized C0 estimate} and \eqref{cm bounds} show that 
	\beq
	\left\|\mu\frac{c_{m+1}}{c_{m}}\hat h_{m}\right\|_{C^0(M,\om_g,h_0)}\leq C_4
	\eeq 
	for a constant $C_4=C_4(\om_g,h_0,\lambda_E,\mu)$. 
	By a similar argument as in the  proof of
	\cite[Theorem~2.5]{XYY26+}, we deduce from \eqref{K bounds} and \eqref{normalized C0 estimate} that 
	\beq \|\hat h_m\|_{C^1(M,\omega_g, h_0)} \leq C_5\eeq 
	where $C_5=C_5(\om_g,h_0,\lambda_E,\mu)$.
Applying the Sobolev embedding theorem and the standard Schauder estimate to the equation
\eqref{normalized iteration equation} yields the uniform bound
\begin{equation}
\|\hat h_m\|_{C^{k}(M,\omega_g, h_0)} \leq \tilde C_k,
\end{equation}
where $\tilde C_k=\tilde C_k(k,\om_g,h_0,\lambda_E,\mu)$ is a constant independent of $m$.
\end{proof}

\vskip 1\baselineskip
\begin{proof}[Proof of Theorem~\ref{main theorem}]
	Fix $\mu>0$ and an initial metric $h_0$ on $E$. By Proposition~\ref{iteration proposition} and Lemma \ref{Fm monotonicity lemma}, there exists a subsequence $\{m_i\}$
	 such that  $ c_{m_i}/c_{m_i-1}\to\alpha$ and
	\beq
	\hat h_{m_i}\to \hat h_\infty,  \ \ \ \hat h_{m_i-1}\to  \hat h_-
	\eeq
 in $C^\infty$,	where $\hat h_\infty$ and $\hat h_-$ are smooth Hermitian metrics on $E$, and $\alpha\in (0,1]$.
	Passing to the limit in \eqref{normalized iteration equation} gives
	\begin{equation}
	S^{\hat h_\infty}=(\lambda_E-\mu)\hat h_\infty+\mu\alpha \hat h_-.
	\label{cluster pair equation}
	\end{equation}
	Define $\tilde h=\alpha \hat h_-$.  Since $K^{\tilde h}=K^{\hat h_-}$,  $\cE(\tilde h)=\cE(\hat h_-)$ where $\cE$  denotes the Hermitian-Yang-Mills energy functional introduced in the proof of Lemma~\ref{adjacent comparison lemma}. Moreover, for every $m\geq 0$,
	\[
	\cE(\hat h_m)=\cE(c_mh_m)=\cE(h_m)=\cE_m.
	\]
	The energy calculation in the proof of Lemma~\ref{adjacent comparison lemma}
	shows that $\{\cE_m\}$ is non-increasing and non-negative, so it has a limit
	$L$.  In particular,
	\beq
	\cE(\hat h_\infty)=\lim_{i\to\infty}\cE_{m_i}=L,
	\qtq{and}
	\cE(\tilde h)=\cE(\hat h_-)=\lim_{i\to\infty}\cE_{m_i-1}=L.\label{cluster energy equality}
	\eeq
	If we set 
	$Q=\tilde h\cdot \hat h_\infty^{-1}$, then \eqref{cluster pair equation} is equivalent to
	\beq K^{\hat h_\infty}-\lambda_E\mathrm{Id}_E=\mu \left(Q-\mathrm{Id}_E\right).\label{keyid}\eeq 
 Since $\cE(\hat h_\infty)=\cE(\tilde h)$,  by the argument in the proof of Lemma~\ref{adjacent comparison lemma} , we obtain 
	\beq \p^{\tilde h} Q=0.
	\eeq 
	Since $Q$ is self-adjoint with respect to
	$\tilde h$,  \beq \bp Q=0.\eeq  
	Since $E$ is stable, $H^0(M, E^*\ts E)\cong \C$. This implies
	\beq
Q=\tilde h\cdot \hat h_\infty^{-1}=a\id_E
	\eeq 
	for some constant $a>0$.  Tracing \eqref{keyid} and integrating over $M$ gives $a = 1$ by the definition of $\lambda_E$.
	Therefore \beq \tilde h=\hat h_\infty, \eeq and \eqref{cluster pair equation} reduces to the Hermitian-Einstein equation
	\beq
	S^{\hat h_\infty}=\lambda_E \hat h_\infty.
	\eeq 
	
	Moreover, suppose that  $\hat h_\infty'$ and $\hat h_\infty''$  are two cluster points of
	$\{\hat h_m\}$.  The preceding paragraph shows that both of them are
	Hermitian-Einstein metrics on $E$. Hence, there exists a constant $b>0$ such that  \beq \hat h_\infty''=b \hat h_\infty'.\eeq  Passing to the limit in \eqref{normalized determinant mean}, we obtain
	\beq 
	0=\int_M\log
	\frac{\det\left(\hat h_\infty''\cdot h_0^{-1}\right)}
	{\det\left(\hat h_\infty'\cdot h_0^{-1}\right)}
	\omega_g^n
	=r\log b\int_M\omega_g^n.
	\eeq 
	Hence $b=1$, and so $\{\hat h_m\}$  converges
	smoothly to the Hermitian-Einstein metric
	$\hat h_\infty$. Moreover, the  Hermitian-Einstein equation can be written as
	\begin{equation}
	S^{\hat h_\infty}=(\lambda_E-\mu)\hat h_\infty+\mu\hat h_\infty.
	\label{fixed point equation}
	\end{equation}
	Choose a constant $C_1=C_1(h_0,\hat h_\infty)$ such that
	\begin{equation}
	h_0\leq C_1\hat h_\infty.
	\label{reduce to original sequence 1}
	\end{equation}
	We claim that for all $m\geq 0$ 
	\beq 
	h_m\leq C_1\hat h_\infty. \label{reduce to original sequence 2}
	\eeq  
	Suppose this  estimate holds for some $m$.  Then one has
	\beq 
	S^{h_{m+1}}-(\lambda_E-\mu)h_{m+1}
	=\mu h_m
	\leq \mu C_1\hat h_\infty
	=S^{C_1\hat h_\infty}-(\lambda_E-\mu)C_1\hat h_\infty .
	\eeq 
	By the comparison theorem \cite[Theorem~2.3]{WYY26b+}, we obtain
	\begin{equation}
	h_{m+1}\leq C_1\hat h_\infty.
	\end{equation}
	Together with  \eqref{normalized C0 estimate}, it follows that 
	\beq
	\tilde C_0^{-1}h_0\leq \hat h_m= c_mh_m\leq c_mC_1\hat h_\infty.
	\eeq
	This implies that 
	\beq 
	c_m\geq C_2
	\eeq 
	for some constant $C_2=C_2(\omega_g, h_0, \lambda_E, \mu)>0$. Since the sequence $\{c_m\}$ is non-increasing (Lemma~\ref{Fm monotonicity lemma}), there exists a constant $c_\infty$ such that
	\begin{equation}
	\lim_{m\to\infty}c_m=c_\infty>0.
	\label{c limit}
	\end{equation}
	Therefore, if we set $h_\infty= c_\infty^{-1}\hat h_\infty$, then
	$	h_m=c_m^{-1}\hat h_m\> h_\infty$ in $C^\infty$.   In particular,  the limit $h_\infty$ satisfies the Hermitian-Einstein equation
	\beq
	\Lambda_{\omega_g}\left(\sqrt{-1}R^{ h_\infty}\right)
	=\lambda_E h_\infty .
	\eeq
	This completes the proof.
\end{proof}

\end{document}